\documentclass{article}
\usepackage{amssymb}
\usepackage{amsfonts,amsmath,epsfig}
\usepackage{mitpress}
\usepackage{graphicx}

\begin{document}

\begin{center}
{\Large Remarks regarding some Algebras of Logic}

\begin{equation*}
\end{equation*}

Mariana Floricica C\u{a}lin, Cristina Flaut, Dana Piciu
\end{center}

\begin{equation*}
\end{equation*}%
\textbf{Abstract.} {\small Algebras of Logic deal with some algebraic
structures, often bounded lattices, considered as models of certain logics,
including logic as a domain of order theory. There are well known their
importance and applications in social life to advance useful concepts, as
for example computer algebra. }

{\small In this paper we reffer in specially to BL-algebras and we present
properties of finite rings or rings with a finite number of ideals in their
connections with BL-rings.}



08C05, 06F35.

\begin{equation*}
\end{equation*}

\bigskip \textbf{0. Introduction}%
\begin{equation*}
\end{equation*}

Algebras of Logic are explicit algebraic systems that show the basic
mathematical structure of Logic. These algebras deal with some algebraic
structures, often bounded lattices, considered as models of certain logics,
including logic as a domain of order theory. There are well known their
importance and applications in social life to advance useful concepts, as
for example computer algebra, also called symbolic computation or algebraic
computation, which it is a scientific area that refers to the study and
development of algorithms or/and software used for computing mathematical
expressions and manipulating other mathematical objects, as well as in the
study of truth values from social analyses. In the following, we will
exemplify some of the applications of Algebras of Logic.

Basic Logic (BL, for short) was introduced by Hajek in [H; 98] to formalize
the many-valued semantics induced on the real interval $[0,1]$ by a
continuous t-norm. Basic Logic generalizes the three most used logics in the
theory of fuzzy sets: \L ukasiewicz logic, G\H{o}del logic and Product
logic. BL-algebras are Lindenbaum-Tarski algebras for Basic Logic.

In this paper we refer especially to BL-algebras and BL-rings.

BL-rings are commutative unitary rings whose lattice of ideals can be
equipped with a structure of BL-algebra. We present properties of finite
rings or rings with a finite number of ideals and their connections with
BL-rings. We give some examples in this purpose and we present in
Conclusions directions of research to find classifications of such rings.

In this paper, all considered rings are commutative and unitary rings.

\begin{equation*}
\end{equation*}

\bigskip \textbf{1. Preliminaries}

\begin{equation*}
\end{equation*}

\bigskip Let $R$ be a commutative unitary ring. The set $Id\left( R\right) $
denotes the set of all ideals of the ring $R$. Let $I,J\in Id\left( R\right) 
$. The following sets are also ideals in $R:$

\begin{equation*}
I+J=<I\cup J>=\{i+j,i\in I,j\in J\}\text{, }
\end{equation*}%
\begin{equation*}
I\otimes J=\{\underset{k=1}{\overset{n}{\sum }}i_{k}j_{k},\text{ }i_{k}\in
I,j_{k}\in J\}\text{, }
\end{equation*}%
\begin{equation*}
\left( I:J\right) =\{x\in R,x\cdot J\subseteq I\}\text{,}
\end{equation*}%
\begin{equation*}
Ann\left( I\right) =\left( \mathbf{0}:I\right) \text{, where }\mathbf{0}=<0>,%
\text{ }
\end{equation*}%
\textbf{\ }called \textit{sum}, \textit{product}, \textit{quotient} and 
\textit{annihilator} of the ideal $I$.\medskip

\textbf{Definition 1.1.} ([AM; 69]) We consider\ $I,J\in Id\left( R\right) $%
. The ideals $I$ and $J$ are called \textit{coprime} if $I+J=R$, that means
there are $i\in I$, $j\in J$ such that $i+j=1$.\textbf{\medskip }

\textbf{Remark 1.2.} If $R_{1}$ and $R_{2}$ are commutaive rings and $%
I\,_{1} $ is an ideal in $R_{1}$ and $I_{2}$ is an ideal in $R_{2},$ then $%
I_{1}\times I_{2}$ is an ideal in the ring $R=R_{1}\times R_{2}$. Moreover,
each ideal in the ring $R$ is on the form $I\times J$, where $I$ is an ideal
in $R_{1}$ and $J$ is an ideal in $R_{2}$.\medskip

\textbf{Proposition 1.3.} ([AM; 69], The Chinese Remainder Theorem for
Rings) \textit{If} $R$ \textit{is a commutative ring and} $I,J$ \textit{two
ideals in }$R$ \textit{such that} $I+J=R$\textit{, we have} $I\otimes
J=I\cap J$. \textit{Then}%
\begin{equation*}
\frac{R}{I\otimes J}=\frac{R}{I\cap J}\simeq \frac{R}{I}\times \frac{R}{J}.
\end{equation*}%
$\Box \medskip $

\textbf{Proposition 1.4.} \textit{1) If }$R_{i},i\in \{1,2,...,s\}$ \textit{%
are commutative unitary rings and} $I_{i}$ \textit{is an ideal in the ring} $%
R_{i}$\textit{, if we consider} $R=R_{1}\times ...\times R_{s}$ \textit{and} 
$I=I_{1}\times ...\times I_{s}$\textit{\ an ideal in the ring }$R,$\textit{\
then} 
\begin{equation*}
\frac{R}{I}\simeq \frac{R_{1}}{I_{1}}\times ...\times \frac{R_{s}}{I_{s}}%
\text{. }
\end{equation*}

2) \textit{For a commutative ring} $R$\textit{, we have} $R\left[ X\right]
\times R\left[ X\right] \simeq (R\times R)\left[ X\right] $.

\textit{3) If} $\left( m,n\right) =1$\textit{, then} $\mathbb{Z}_{m}\left[ X%
\right] \times \mathbb{Z}_{n}\left[ X\right] \simeq (\mathbb{Z}_{m}\times 
\mathbb{Z}_{n})\left[ X\right] =\mathbb{Z}_{mn}\left[ X\right] $.

\textit{4)} 
\begin{equation*}
\frac{\mathbb{Z}\left[ X\right] }{\left( n\right) }\simeq (\mathbb{Z}/\left(
n\right) )[X].
\end{equation*}

\textit{5)} \textit{If} $I=I_{1}\times I_{2}$ \textit{is an ideal in the ring%
} $R\left[ X\right] \times R\left[ X\right] $\textit{, then} 
\begin{equation*}
(R\left[ X\right] \times R\left[ X\right] )/I\simeq \frac{R\left[ X\right] }{%
I_{1}}\times \frac{R\left[ X\right] }{I_{2}}.
\end{equation*}

\textbf{Proof.} 1) It is clear from the Fundamenthal Theorem of Isomorphism,
considering the surjective morphism $\phi :$ $\frac{R}{I}\rightarrow \frac{%
R_{1}}{I_{1}}\times ...\times \frac{R_{s}}{I_{s}},\phi \left(
x_{1},...,x_{s}\right) =\left( \widehat{x_{1}},...,\widehat{x_{s}}\right) $.

2) Let $f:R\left[ X\right] \times R\left[ X\right] \simeq (R\times R)\left[ X%
\right] ,f\left( p,q\right) =Q[X],$ where $p=p_{0}+p_{1}X+...+p_{m}X^{m}$, $%
q=q_{0}+q_{1}X+...+q_{m}X^{m}$, with $q_{m}$ or $p_{m}\,\ $possible zero and 
$Q[X]=\left( p_{0},q_{0}\right) +\left( p_{1},q_{1}\right) X+...+\left(
p_{m},q_{m}\right) X^{m}$. It is clear that $f$ is a bijective morphism.

3) Indeed, let $p\left( X\right) \in \mathbb{Z}\left[ X\right] $ be a
polynomial and $f:$ $\mathbb{Z}_{mn}\left[ X\right] \rightarrow \mathbb{Z}%
_{m}\left[ X\right] \times \mathbb{Z}_{n}\left[ X\right] ,$

$f\left( \overline{p\left( X\right) }\right) =\left( \widehat{p_{m}\left(
X\right) },\widetilde{p_{n}\left( X\right) }\right) $, where $\overline{%
p\left( X\right) }$ is the polynomial $p\left( X\right) $ with the
coefficients reduced mod $mn$, $\widehat{p_{m}\left( X\right) }~$\ is the
polynomial $p\left( X\right) $ with the coefficients reduced mod $m$ and $%
\widetilde{p_{n}\left( X\right) }$ \ is the polynomial $p\left( X\right) $
with the coefficients reduced mod $n$. We have, $f\left( \overline{p_{1}}+%
\overline{p_{2}}\right) =f\left( \overline{p_{1}}\right) +f\left( \overline{%
p_{2}}\right) $ and $f\left( \overline{p_{1}}\overline{p_{2}}\right)
=f\left( \overline{p_{1}}\right) f\left( \overline{p_{2}}\right) $.

4) We have $\pi :\mathbb{Z\rightarrow Z}/\left( n\right) $, canonical
projection and $i:\mathbb{Z}/\left( n\right) \rightarrow \mathbb{Z}/\left(
n\right) \left[ X\right] $ canonical injection. Therefore, from universal
propertie of polinomial rings, we have the morphism ring $f:\mathbb{Z}\left[
X\right] \rightarrow \mathbb{Z}/\left( n\right) \left[ X\right] $ with 
\textit{ker}$f=n\mathbb{Z}\left[ X\right] =\left( n\right) =\{$ the set of
polynomials with coefficients multiple of $n\}$. It is clear that $f$ is
surjective. Therefore, from Isomorphism Fundamental Theorem, we have that $%
\frac{\mathbb{Z}\left[ X\right] }{\left( n\right) }\simeq (\mathbb{Z}/\left(
n\right) )[X]$.

5) It is clear from 1). $\Box \medskip $

\textbf{Definition 1.5}. ([WD; 39]) \ A \textit{(commutative) residuated
lattice }\ is an algebra $(L,\wedge ,\vee ,\odot ,\rightarrow ,0,1)$ such
that:

(1) $\ (L,\wedge ,\vee ,0,1)$ is a bounded lattice;

(2) $\ \ (L,\odot ,1)$ is a commutative ordered monoid;

(3) $\ z\leq x\rightarrow y$\ iff $x\odot z\leq y,$\ for all $x,y,z\in L.$

The property (3) is called\emph{\ }\textit{residuation}, where $\leq $ is
the partial order of the lattice $\ (L,\wedge ,\vee ,0,1).$

In a residuated lattice is defined an additional operation: for $x\in L,$ we
denote $x^{\ast }=x\rightarrow 0.$

If we preserve these notations, for a commutative and unitary ring we have
that $(Id(R),\cap ,+,\otimes \rightarrow ,0=\{0\},1=R)$ is a residuated
lattice in which the order relation is $\subseteq $, $I\rightarrow J=(J:I)$
and $I\odot J=I\otimes J,$ for every $I,J\in Id(R)$, see [TT;22]

In a residuated lattice $(L,\wedge ,\vee ,\odot ,\rightarrow ,0,1)$ we
consider the following identities:

\begin{equation*}
(prel)\qquad (x\rightarrow y)\vee (y\rightarrow x)=1\qquad \text{ (\textit{%
prelinearity)}};
\end{equation*}%
\begin{equation*}
(div)\qquad x\odot (x\rightarrow y)=x\wedge y\qquad \text{ (\textit{%
divisibility)}}.
\end{equation*}%
\textbf{Definition 1.6. }([T; 99]) 1) A residuated lattice $L$ is called 
\textit{a BL-algebra\ }if $L$ verifies $(prel)$ and $(div)$ conditions.

2) A \textit{BL-chain}\emph{\ }is a totally ordered BL-algebra, i.e., a
BL-algebra such that its lattice order is total.\medskip \medskip

\textbf{Definition 1.7. }([CHA; 58])\textbf{\ }An \textit{MV-algebra }is an
algebra $\left( L,\oplus ,^{\ast },0\right) $ satisfying the following
axioms:

(1) $\left( L,\oplus ,0\right) $ \ is an abelian monoid;

(2) $(x^{\ast })^{\ast }=x;$

(3) $x\oplus 0^{\ast }=0^{\ast };$

(4) $\left( x^{\ast }\oplus y\right) ^{\ast }\oplus y=$ $\left( y^{\ast
}\oplus x\right) ^{\ast }\oplus x$, for all $x,y\in L.\medskip $

\textbf{Remark 1.8.} If in a BL- algebra $L$ we have $x^{\ast \ast }=x,$ for
every $x\in L$, and, for $x,y\in L,$ we denote 
\begin{equation*}
x\oplus y=(x^{\ast }\odot y^{\ast })^{\ast },
\end{equation*}%
then we obtain an MV-algebra structure $(L,\oplus ,^{\ast },0).$ Conversely,
if $(L,\oplus ,^{\ast },0)$ is an MV\textit{-}algebra, then $(L,\wedge ,\vee
,\odot ,\rightarrow ,0,1)$ becomes a BL-algebra, in which for $x,y\in L$ we
have$:$ 
\begin{equation*}
x\odot y=(x^{\ast }\oplus y^{\ast })^{\ast },
\end{equation*}%
\begin{equation*}
x\rightarrow y=x^{\ast }\oplus y,1=0^{\ast },
\end{equation*}%
\begin{equation*}
x\vee y=(x\rightarrow y)\rightarrow y=(y\rightarrow x)\rightarrow x\text{
and }x\wedge y=(x^{\ast }\vee y^{\ast })^{\ast },
\end{equation*}%
see [T; 99].

We recall that in [NL; 03], Di Nolla and Lettieri analyze the structure of
finite BL-algebras. They introduced the concept of BL-comets, a class of
finite BL-algebras which can be seen as a generalization of finite
BL-chains. Using BL-comets, any finite BL-algebra can be represent as a
direct product of BL-comets.

\textbf{Definition 1.9.} 1) ([NL; 03], Definition 21) Let $L$ be a
BL-algebra. An element $x\in L$ is called \textit{idempotent} if $x\odot x=x$%
.

2) Let $L$ be a finite BL-algebra and $\mathcal{I}\left( L\right) $ the set
of idempotent elements in $L$. For $x\in \mathcal{I}\left( L\right) $, we
denote $\mathcal{C}\left( x\right) =\{y\in \mathcal{I}\left( L\right) ~$\
such that $x$ and $y$ are comparable$\}$. We defined the set $\mathcal{D}%
\left( L\right) \subseteq \mathcal{I}\left( L\right) $ as follows:

$x\in \mathcal{D}\left( L\right) $ if and only if

i) $\mathcal{C}\left( x\right) =\mathcal{I}\left( L\right) ;$

ii) The set $\{y\in \mathcal{I}\left( L\right) ,$ $y\leq x\}$ is a chain.

It is clear that $\mathcal{D}\left( L\right) \neq \emptyset $ since $0\in 
\mathcal{D}\left( L\right) .$

A finite BL-algebra $L$ is called a \textit{BL-comet} if \textit{max}$%
\mathcal{D}\left( L\right) \neq 0.\medskip $

We recall that in a BL-comet $L,$ the element \textit{max}$\mathcal{D}\left(
L\right) $ is called the pivot of $L$ and it is denoted by $%
pivot(L).\medskip $

\textbf{Proposition 1.10. }([NL; 03], Proposition 26) \textit{Let} $L$ 
\textit{be a finite BL-algebra. The following assertions are equivalent:}

(i) $L$ is a \textit{BL-comet and }$pivot(L)=1;$

(ii)\textit{\ }$L$ is a BL-chain.

Using the characterization of Boolean elements in BL-algebras, we establish
the connections between BL-comets and BL(MV)-chains.\medskip

\textbf{Proposition 1.11.} \textit{Let} $L$ \textit{be a BL-comet. Then }$L$
is a BL-chain iff $pivot(L)^{\ast \ast }=pivot(L).$

\textbf{Proof.} If $L$ is a BL-chain, using Proposition 1.10, $pivot(L)=1$ \
and \ obviously, $1^{\ast \ast }=1.$

Conversely, if we suppose that $pivot(L)^{\ast \ast }=pivot(L),$ using the
characterizations of Boolean elements in BL-algebras, see [P; 07,
Proposition 3.3], it is clear that $pivot(L)$ is a boolean element in $L,$
so, $pivot(L)\vee pivot(L)^{\ast }=1.$ We deduce that $pivot(L)^{\ast }$ \
is also a boolean element, so, an idempotent element in $L.$ Then $%
pivot(L)^{\ast }\in C(pivot(x))$. We conclude that $pivot(L)$ and $%
pivot(L)^{\ast }$ are comparable.

If $pivot(L)\leq pivot(L)^{\ast },$ then $1=pivot(L)\vee pivot(L)^{\ast
}=pivot(L)^{\ast },$ so, $pivot(L)=0,$ a contradiction.

If $pivot(L)^{\ast }\leq pivot(L),$ then $1=pivot(L)\vee pivot(L)^{\ast
}=pivot(L).$ Using Proposition 1.10, $L$ is a BL-chain. $\Box \medskip $

From Proposition 1.11 we deduce the following result:

\textbf{Proposition 1.12.} \textit{Let} $L$ \textit{be a finite MV-algebra.
The following assertions are equivalent:}

(i) $L$ is a \textit{BL-comet;}

(ii)\textit{\ }$L$ is an MV-chain.

\begin{equation*}
\end{equation*}

\textbf{2.} \textbf{Some remarks regarding commutative unitary rings}%
\begin{equation*}
\end{equation*}

\textbf{Definition 2.1. }([AM; 69]) Let $R$ be a commutative unitary ring.

1) The ideal $M$ of the ring $R$ is \textit{maximal} if it is maximal, with
respect of the set inclusion, amongst all proper ideals of the ring $R$.
That means, there are no other ideals different from $R$ contained $M$. The
ideal $J$ of the ring $R$ is a \textit{minimal ideal} if it is a nonzero
ideal which contains no other nonzero ideals.

2) A commutative \textit{local ring} $R$ is a ring with a unique maximal
ideal.

3) Let $P$\thinspace $\not=R$ be an ideal in the ring $R$. Let $a,b\in R$
such that $ab\in P$. If we have \thinspace $a\in P$ or $b\in P$, therefore $%
P $ is called a \textit{prime} ideal of $R$.\medskip\ 

\textbf{Remark 2.2.} Let $R$ be a commutative unitary ring.

1) The ideal $M$ of the ring $R$ is maximal if and only if $R/M$ is a field;

2) The ideal $P$ of the ring $R$ is prime if and only if $R/P$ is an
integral domain.

From here, we have that a maximal ideal it is a prime ideal.

For other details and properties, the reader is referred to [AM; 69].\medskip

\textbf{Definition 2.3.} i) [HLN; 18] A commutative ring $R$ is called a 
\textit{multiplication ring} if for every ideals $I_{1},I_{2}$ of $R$ with $%
I_{1}\subseteq I_{2}$, there is an ideal $I_{3}$ of $R$ such that $%
I_{1}=I_{2}\otimes I_{3}$.

ii) A commutative ring $R$ is called \textit{Noetherian ring} if the
condition of ascending chain is satisfied, that means every increasing
sequence of ideals $I_{1}\subseteq I_{2}\subseteq ...\subseteq
I_{r}\subseteq ...$ is stationary, that means there is $q$ such that $%
I_{q}=I_{q+1}=...$.

iii) A commutative ring $R$ is called \textit{Artinian ring} if the
condition of descending chain is satisfied, that means every decreasing
sequence of ideals $I_{1}\supseteq I_{2}\supseteq ...\supseteq
I_{r}\supseteq ...$ is stationary, that means there is $q$ such that $%
I_{q}=I_{q+1}=...$.\medskip

\textbf{Remark 2.4.}

i) ([AB; 19], Lemma 3.5) Let $R=\underset{j\in J}{\prod }R_{i}$ be a direct
product of rings. $R$ is a multiplication ring if and only if $R_{j}$ is a
multiplication ring for all $j\in J$.

ii) ([AB; 19], Lemma 3.6) Let $R$ be a multiplication ring and $I$ be an
ideal of $R$. Therefore, \thinspace the quotient ring $R/I$ is a
multiplication ring.\medskip

\textbf{Remark 2.5.} ( [A; 76]. Corollary 6.1) Let $R$ be a ring. The
following conditions are equivalent:

i) $R[X]$ is a multiplication ring;

ii) $R$ is a finite direct product of fields.\medskip

\textbf{Remark 2.6.} ([AF; 92] and [AM; 69])

1) If $R$ is a Noetherian ring, therefore the polynomial ring $R\left[ X%
\right] $ is Noetherian and the quotient ring $R/I$ is also a Noetherian
ring, for $I~$an ideal of $R$.

2) Any field and any principal ideal ring is a Noetherian ring.

3) Every ideal of the Noetherian ring $R$ is finitely generated.

4) An integral domain $R$ is Artinian ring if and only if $R$ is a field.

5) The ring $K\left[ X\right] /\left( X^{t}\right) $ is Artinian ring, for $%
K $ a field and $t$ a positive integer.

6) A commutative Noetherian ring $R$ is Artinian if and only if $R$ is a
product of local rings.

7) In an Artinian ring every prime ideal is maximal.

8) An Artinian ring is a finite direct product of Artinian local
rings.\medskip

\textbf{Proposition 2.7.} \textit{Let} $R$ \textit{be a commutative unitary
ring which has exactly three ideals }$\{0\},I,R$. \textit{Therefore, we have}
$I^{2}=\{0\}$.\medskip

\textbf{Proof.} It is clear that $I$ is the only minimal ideal in $R$. Since 
$I^{2}\subseteq I$, we have $I^{2}=\{0\}$ or $I^{2}=I$. In the same time, $I$
is the only maximal ideal of $R$ and each element from $I$ is a non-unit.
Suppose that $I^{2}=I$. Let $x\in I$ such that $x\neq 0$. It is clear that $%
I=<x>$ and $x^{2}\neq 0$. We consider the ideal $A=\{z\in R~/~xz=0\}\,$.We
prove that $A\neq \{0\}$. Since $x^{2}\neq 0$, therefore there is an element 
$y\in R$ such that $x^{2}y=x$. It results \thinspace \thinspace $x\left(
xy-1\right) =0,$ with $xy-1\neq 0,$ since $x$ is a non unit element in $R$.
Then $xy-1\in A$ and $A\neq \{0\}$. We have that $A\neq R$ since $x\notin A$
( we have $x^{2}\neq 0$). From here, we obtain that $A=I,$ a contradiction,
therefore $I^{2}=\{0\}$.$\Box \medskip $

\textbf{Proposition 2.8.} \textit{Each nonzero element in a finite
commutative unitary ring }$R$\textit{\ is a unit or a zero divisor.\medskip }

\textbf{Proof.} Let $R=\{x_{1},x_{2},...,x_{n}\}$ be a ring with $n$
elements and $a\in R$ a non unit element. Therefore, the following elements $%
ax_{1},ax_{2},...,ax_{n}$ are all different from $1$. It results that we can
have maximum $n-1$ distict values, then there are equal at least two of
them. \ We otain that $ax_{i}=ax_{j}$, for $i\neq j$, therefore $a\left(
x_{i}-x_{j}\right) =0$. Since $x_{i}\neq x_{j}$, it results that $a$ is a
zero divisor.$~\Box \medskip $

\textbf{Remark 2.9.} (see [FK; 12])

1. In a commutative ring $R$, the set of non-unit elements is an ideal if
and only if the ring $R$ is local.

2. If $R$ is a commutative Noetherian ring and $I$ is an ideal in $R$
consists of zero-divisors, then its annihilator is a non zero ideal.\medskip

\textbf{Proposition 2.10.} \textit{Let} $R$ \textit{be an integral domain
with a finite number of ideals. Then} $R$ \textit{is a field.\medskip\ }

\textbf{Proof.} Let $a\in R$ be a nonzero element. We have the ideals $%
<a>,<a^{2}>,...,<a^{k}>,...$Since we have a finite number of ideals, it
results that we can find two equal ideals in this chain: $%
<a^{k_{1}}>=<a^{k_{2}}>$. Therefore, there is an element $r\in R$ such that $%
a^{k_{1}}r=a^{k_{2}}$. Assuming that $k_{1}\geq k_{2}$, it results that $%
a^{k_{2}}\left( a^{k_{1}-h_{2}}r-1\right) =0$. Since $a\not=0$, we have that 
$a^{k_{1}-h_{2}}r-1=0$, therefore $a$ is an invertible element. $\Box
\medskip $

\textbf{Proposition 2.11.} \textit{Let} $R$ \textit{be a commutative and
unitary ring with a finite number of ideals. Therefore, an ideal} $I$ 
\textit{is maximal if and only if it is a prime ideal.\medskip }

\textbf{Proof.} Let $I$ be a maximal ideal. Therefore, $R/I$ $\ $is a field,
that means an integral domain, therefore $I$ is a prime ideal. Conversely,
if $I$ is a prime ideal, therefore $R/I$ is an integral domain with a finite
number of ideals. From the above proposition, $R/I$ is a field, then $I$ is
a maximal ideal. $\Box \medskip $

\textbf{Proposition 2.12.} \textit{Let} $R$ \textit{be a commutative and
unitary ring with a finite number of ideals. Let} $n_{m}\left( R\right) $ 
\textit{be the number of maximal ideals in} $R$\textit{,} $n_{p}\left(
R\right) $ \textit{be the number of prime ideals in} $R$ \textit{and} $%
n_{I}\left( R\right) $ \textit{be the number of all ideals in} $R$\textit{.
Therefore,} $n_{m}\left( R\right) =n_{p}\left( R\right) =\alpha $ \textit{and%
} $n_{I}\left( R\right) =\underset{j=1}{\overset{\alpha }{\prod }}\beta
_{j},\beta _{j}$ \textit{positive integers, }$\beta _{j}\geq 2$.

\textbf{Proof.} Relation $n_{m}\left( R\right) =n_{p}\left( R\right) =\alpha 
$ arises from the above proposition. It is clear that $R$ is an Artinian
ring, therefore it is a product of local nonzero rings, $R=R_{1}\times
...\times R_{\alpha }$. Therefore, each ideal $J$ in $R$ is of the form $%
J=J_{1}\times ...\times J_{\alpha }$, with $J_{i}$ ideal in $R_{i}$. Since
each ring $R_{i}$ has at least two ideals, we obtain that $n_{I}\left(
R\right) =\underset{j=1}{\overset{\alpha }{\prod }}\beta _{j}$, where $\beta
_{j}\geq 2$ is the number of ideals in $R_{j}$. $\Box \medskip $

\textbf{Remark} \textbf{2.13}. 1) First, we must remark that $n_{I}\left(
R_{j}\right) $ can be any positive integer $\beta _{j}\geq 2\,$. For
example, the ring $K\left[ X\right] /\left( X^{\beta _{j}}\right) $ has $%
\beta _{j}$ ideals, $\beta _{j}\geq 2$.

2) If $n_{I}\left( R\right) $ is finite, from Proposition 2.11, we have only
the following two possibilities:

-$R$ is an integral domain, therefore it is a field and, in this case, we
have $n_{m}\left( R\right) =n_{p}\left( R\right) =1$ and $n_{I}\left(
R\right) =2~$or

- $R$ is not an integral domain and $n_{m}\left( R\right) =n_{p}\left(
R\right) \geq 1$ and $n_{I}\left( R\right) >2$.

3) From the above proposition it is clear that there are not commutative
unitary rings such that $\left( n_{m}\left( R\right) ,n_{p}\left( R\right)
,n_{I}\left( R\right) \right) =\left( 3,3,5\right) $ or $\left( n_{m}\left(
R\right) ,n_{p}\left( R\right) ,n_{I}\left( R\right) \right) =\left(
2,2,5\right) $, since $5$ is a prime number. To find such an example, we
must search in non-commutative rings. Therefore we have examples only in the
case $\left( n_{m}\left( R\right) ,n_{p}\left( R\right) ,n_{I}\left(
R\right) \right) =\left( 1,1,5\right) $. The same situation for $n_{I}\left(
R\right) =7$. We have only the case $\left( n_{m}\left( R\right)
,n_{p}\left( R\right) ,n_{I}\left( R\right) \right) =\left( 1,1,7\right) $
But, we can find examples of commutative unitary rings such that $\left(
n_{m}\left( R\right) ,n_{p}\left( R\right) ,n_{I}\left( R\right) \right)
=\left( 2,2,4\right) $ or commutative unitary rings such that $\left(
n_{m}\left( R\right) ,n_{p}\left( R\right) ,n_{I}\left( R\right) \right)
=\left( 1,1,4\right) $. For example, for the first case, we have the ring $%
\mathbb{Z}_{2}\times \mathbb{Z}_{2}$ or $\mathbb{Z}_{6}$ and for the second,
we have the ring $\mathbb{Z}_{8}.$

\begin{equation*}
\end{equation*}

\textbf{3.} \textbf{Remarks regarding} \textbf{BL-rings}%
\begin{equation*}
\end{equation*}

Let $R$ be a commutative unitary ring and $Id\left( R\right) $ be the set of
all ideals of the ring $R$. As we remarked above, we have that $(Id(R),\cap
,+,\otimes \rightarrow ,0=\{0\},1=R)$ is a residuated lattice in which the
order relation is $\subseteq $ and $I\rightarrow J=(J:I),$ for every $I,J\in
Id(R)$. A commutative ring is a \textit{BL-ring} if and only if $Id(R)$ is a
BL-algebra (see [HLN; 18], Corollary 2.3.). A Noetherian multiplicative ring
is an example of BL-ring (see [HLN; 18], Example 2.4, 2.).

If we consider $R_{1}$ and $R_{2}$ two BL-rings and their lattices of ideals 
$Id\left( R_{1}\right) $ and $Id\left( R_{2}\right) $, therefore the lattice 
$Id\left( R_{1}\times R_{2}\right) =Id\left( R_{1}\right) \times Id\left(
R_{2}\right) $ is a BL-algebra, therefore $R_{1}\times R_{2}$ is a BL-ring.
Indeed, for $Z,W\in Id\left( R_{1}\times R_{2}\right) ,Z=\left(
I_{1},J_{1}\right) $ and $W=\left( I_{2},J_{2}\right) $, we define 
\begin{equation*}
Z\rightarrow W=\left( I_{1}\rightarrow I_{2},J_{1}\rightarrow J_{2}\right) ,
\end{equation*}%
and 
\begin{equation*}
Z\otimes W=\left( I_{1}\otimes J_{1},I_{2}\otimes J_{2}\right) .
\end{equation*}%
In this way, we obtain a BL-algebra structure on $Id\left( R_{1}\times
R_{2}\right) $.

\textbf{Proposition 3.1.} \textit{There are no commutative unitary rings} $R$
\textit{with three ideals having} $(Id(R),\cap ,+,\otimes \rightarrow
,0=\{0\},1=R)$ \textit{as a BL-algebra which is not an MV-algebra.\medskip }

\textbf{Proof.} From Proposition 2.7, it results that if $R$ is a
commutative unitary ring which has exactly three ideals $\{0\},I,R$,
therefore we have $I^{2}=\{0\}$. \ From here, we have $I\rightarrow 0=\left(
0:I\right) =I$ and the following implication and multiplication tables are
obtained

\begin{equation}
\begin{tabular}{l|lll}
$\rightarrow $ & $0$ & $I$ & $R$ \\ \hline
$0$ & $R$ & $R$ & $R$ \\ 
$I$ & $I$ & $R$ & $R$ \\ 
$R$ & $0$ & $I$ & $R$%
\end{tabular}%
\ \ \ 
\begin{tabular}{l|lll}
$\otimes $ & $0$ & $I$ & $R$ \\ \hline
$0$ & $0$ & $0$ & $0$ \\ 
$I$ & $0$ & $0$ & $I$ \\ 
$R$ & $0$ & $I$ & $R$%
\end{tabular}%
\ \ ,  \tag{3.1.}
\end{equation}%
which give us an MV-algebra structure. An example of such a ring is $\left( 
\mathbb{Z}_{4},+,\cdot \right) $, with $I=\left( \widehat{2}\right)
,0=\left( \widehat{0}\right) $. $\Box \medskip $

\textbf{Remark 3.2.} From the above, a commutative unitary ring $R$\ with
three ideals always has on the algebra $Id\left( R\right) $\ an MV-algebra
structure. Therefore, to obtain rings with three ideals such that $Id\left(
R\right) $ is a BL-algebra that is not an MV-algebra, we must search among
non-unitary rings or among non-commutative rings or both, or to find such an
algebra as a subset or as a subalgebra of a infinite BL-algebra. The only
BL-algebra of order $3$ which is not an MV-algebra has the following
implication and multiplication tables:%
\begin{equation}
\begin{tabular}{l|lll}
$\rightarrow $ & $0$ & $I$ & $R$ \\ \hline
$0$ & $R$ & $R$ & $R$ \\ 
$I$ & $0$ & $R$ & $R$ \\ 
$R$ & $0$ & $I$ & $R$%
\end{tabular}%
\ \ \ 
\begin{tabular}{l|lll}
$\otimes $ & $0$ & $I$ & $R$ \\ \hline
$0$ & $0$ & $0$ & $0$ \\ 
$I$ & $0$ & $I$ & $I$ \\ 
$R$ & $0$ & $I$ & $R$%
\end{tabular}%
.  \tag{3.2.}
\end{equation}

\textbf{Remark 3.3. }1) We recall that if $K$ is a field and $K\left[ X%
\right] $ is the polynomial ring, then for $f\in K\left[ X\right] $, the
quotient ring $A=K\left[ X\right] /\left( f\right) $ is a BL-ring, see [FP;
22].

2) Let $n=p_{1}p_{2}$, with $p_{1}$ and $p_{2}$ two prime distinct integers.
We have the following isomorphism 
\begin{equation*}
R=\frac{\mathbb{Z}_{p_{1}p_{2}}\left[ X\right] }{\left(
p_{1}X^{2}+p_{2}\right) }\simeq \frac{\mathbb{Z}\left[ X\right] }{\left(
n,~p_{1}X^{2}+p_{2}\right) }\simeq \frac{\mathbb{Z}\left[ X\right] }{\left(
p_{2},~X^{2}\right) }\simeq \frac{\mathbb{Z}_{p_{2}}\left[ X\right] }{\left(
X^{2}\right) }.\medskip
\end{equation*}%
This ring is a BL-ring.\medskip

\textbf{Example 3.4.} In the following, we present some examples of
BL-rings, more precisely, examples of Noetherian multiplicative rings which
are BL-rings.

1) We consider the ring $R=\frac{\mathbb{Z}_{6}\left[ X\right] }{\left(
X^{2}\right) }$. From the above, we have $\frac{\mathbb{Z}_{6}\left[ X\right]
}{\left( X^{2}\right) }\simeq \frac{\mathbb{Z}_{2}\left[ X\right] }{\left(
X^{2}\right) }\times \frac{\mathbb{Z}_{3}\left[ X\right] }{\left(
X^{2}\right) }$. The ring $\frac{\mathbb{Z}_{2}\left[ X\right] }{\left(
X^{2}\right) }$ has three ideals $\left( 0\right) ,\left( X\right) $ and $%
\frac{\mathbb{Z}_{2}\left[ X\right] }{\left( X^{2}\right) }$. The ring $%
\frac{\mathbb{Z}_{3}\left[ X\right] }{\left( X^{2}\right) }$ has also three
ideals $\left( 0\right) ,\left( X\right) $ and $\frac{\mathbb{Z}_{3}\left[ X%
\right] }{\left( X^{2}\right) }$ Therefore, the ring $\frac{\mathbb{Z}_{6}%
\left[ X\right] }{\left( X^{2}\right) }$ is a BL-ring with $9$ ideals.

2) From the above results, we have the following rings isomorphisms 
\begin{equation*}
R_{1}=\frac{\mathbb{Z}_{6}\left[ X\right] }{\left( 2X^{2}+3\right) }\simeq 
\frac{\mathbb{Z}\left[ X\right] }{\left( 6,~2X^{2}+3\right) }\simeq \frac{%
\mathbb{Z}\left[ X\right] }{\left( 3,~X^{2}\right) }\simeq \frac{\mathbb{Z}%
_{3}\left[ X\right] }{\left( X^{2}\right) }
\end{equation*}%
and 
\begin{equation*}
R_{2}=\frac{\mathbb{Z}_{6}\left[ X\right] }{\left( 3X^{2}+2\right) }\simeq 
\frac{\mathbb{Z}\left[ X\right] }{\left( 6,~3X^{2}+2\right) }\simeq \frac{%
\mathbb{Z}\left[ X\right] }{\left( 2,~X^{2}\right) }\simeq \frac{\mathbb{Z}%
_{2}\left[ X\right] }{\left( X^{2}\right) }\text{.}
\end{equation*}

Indeed, since $X^{2}=\left( 2X^{2}+3\right) \left( 3X^{2}+2\right) $, we
have that $\left( 3,~X^{2}\right) \subseteq \left( 6,~2X^{2}+3\right) $.
Conversely it is also true, since $2X^{2}+3=2\cdot X^{2}+3$ and $6=3+3$.
Therefore $\left( 6,~2X^{2}+3\right) =\left( 3,~X^{2}\right) $. Similar, we
have $\left( 6,~3X^{2}+2\right) =\left( 2,~X^{2}\right) $. In the rings $%
R_{1}$ and $R_{2},$ we have $X^{2}=0,$ since $X^{2}=\left( 2X^{2}+3\right)
\left( 3X^{2}+2\right) $. Therefore, in $R_{1}$, we have $3=0$ and in $%
R_{2}\,\ $we have $2=0.$

3) From the above results, we have the following isomorphism of rings: 
\begin{equation*}
\frac{\mathbb{Z}_{6}\left[ X\right] }{\left( 2X^{2}+4\right) }\simeq \frac{%
\mathbb{Z}\left[ X\right] }{\left( 6,~2X^{2}+4\right) }\simeq \frac{\mathbb{Z%
}\left[ X\right] }{\left( 2\right) \left( 3,~X^{2}+2\right) }\simeq
\end{equation*}%
\begin{equation*}
\simeq \frac{\mathbb{Z}\left[ X\right] }{\left( 2\right) }\times \frac{%
\mathbb{Z}\left[ X\right] }{\left( 3,~X^{2}+2\right) }\simeq \mathbb{Z}_{2}%
\left[ X\right] \times \frac{\mathbb{Z}_{3}\left[ X\right] }{\left(
X^{2}+2\right) }.
\end{equation*}%
Therefore, the ring $\frac{\mathbb{Z}_{6}\left[ X\right] }{\left(
2X^{2}+4\right) }$ is an infinite BL-ring.

In the same way, we have the following isomorphisms of rings:

\begin{equation*}
\frac{\mathbb{Z}_{30}\left[ X\right] }{\left( 6X^{2}+24\right) }\simeq \frac{%
\mathbb{Z}\left[ X\right] }{\left( 30,~6X^{2}+24\right) }\simeq \frac{%
\mathbb{Z}\left[ X\right] }{\left( 6\right) \left( 5,~X^{2}+4\right) }\simeq
\end{equation*}%
\begin{equation*}
\simeq \frac{\mathbb{Z}\left[ X\right] }{\left( 6\right) }\times \frac{%
\mathbb{Z}\left[ X\right] }{\left( 5,~X^{2}+4\right) }\simeq \mathbb{Z}_{6}%
\left[ X\right] \times \frac{\mathbb{Z}_{5}\left[ X\right] }{\left(
X^{2}+4\right) },
\end{equation*}%
or 
\begin{equation*}
\frac{\mathbb{Z}_{30}\left[ X\right] }{\left( 6X^{2}\right) }\simeq \frac{%
\mathbb{Z}\left[ X\right] }{\left( 30,~6X^{2}\right) }\simeq \frac{\mathbb{Z}%
\left[ X\right] }{\left( 6\right) \left( 5,~X^{2}\right) }\simeq
\end{equation*}%
\begin{equation*}
\simeq \frac{\mathbb{Z}\left[ X\right] }{\left( 6\right) }\times \frac{%
\mathbb{Z}\left[ X\right] }{\left( 5,~X^{2}\right) }\simeq \mathbb{Z}_{6}%
\left[ X\right] \times \frac{\mathbb{Z}_{5}\left[ X\right] }{\left(
X^{2}\right) }.
\end{equation*}%
Therefore, the rings $\mathbb{Z}_{6}\left[ X\right] \times \frac{\mathbb{Z}%
_{5}\left[ X\right] }{\left( X^{2}+4\right) }$ and $\mathbb{Z}_{6}\left[ X%
\right] \times \frac{\mathbb{Z}_{5}\left[ X\right] }{\left( X^{2}\right) }$
are infinite BL-rings.

4) From the above results, we have the following isomorphism of rings: 
\begin{equation*}
\frac{\mathbb{Z}_{6}\left[ X\right] }{\left( X^{2}-X\right) }\simeq \frac{%
\mathbb{Z}_{6}\left[ X\right] }{X\left( X-1\right) }\simeq \frac{\mathbb{Z}%
_{6}\left[ X\right] }{\left( X\right) }\times \frac{\mathbb{Z}_{6}\left[ X%
\right] }{\left( X-1\right) }\simeq
\end{equation*}%
\begin{equation*}
\simeq \mathbb{Z}_{6}\times \mathbb{Z}_{6}.
\end{equation*}

\textbf{Example 3.5.} Let $R=\mathbb{Z}_{6}\left[ X\right] \times \frac{%
\mathbb{Z}_{5}\left[ X\right] }{\left( X^{2}\right) }$, be a BL-ring from
the above examples. Therefore $Id\left( R\right) $ is an infinite
BL-algebra. In $R,$ we consider the following ideals $\{0\},I=\left(
3\right) \times R$ and $R,$ where $\left( 3\right) =\{f\in \mathbb{Z}_{6}%
\left[ X\right] ,f=a_{n}X^{n}+...+a_{1}X+a_{0},$ with $a_{0}=0$ or $%
a_{0}=3\} $. We have that $I^{2}=I$ and $Ann\left( I\right) =\{0\}$.
Therefore, $B=\{\{0\},I,R\}$ is a finite BL-subalgebra of the algebra $%
Id\left( R\right) $, having the implication and multiplication tables given
by the relation $\left( 3.2.\right) $.\medskip

\textbf{Remark 3.6.} Let $R$ be a commutative unitary ring which has exactly
four ideals $\{0\},I,J,R$. Therefore, the lattice $Id\left( R\right) $ can
be of the form $A$ or $B$ from Figure 1. 
\begin{figure}[tbph]
\centerline{\includegraphics[width=2.5in, height=1.3in]{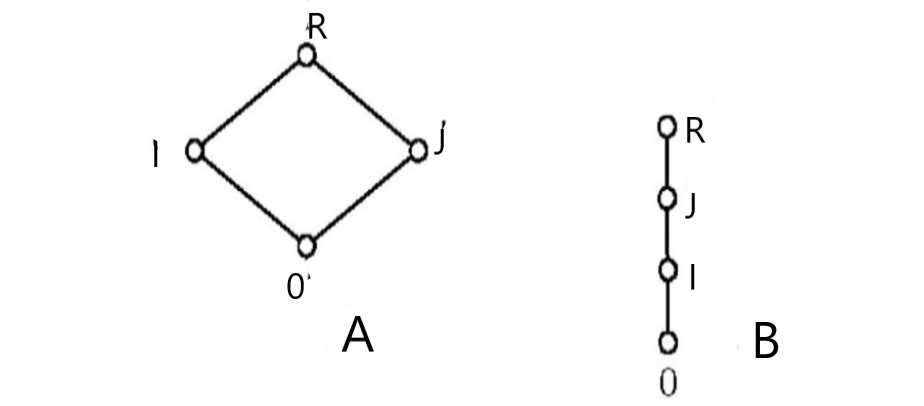}}
\caption{Latices with four elements.}
\label{fig}
\end{figure}
\textbf{Case i}. $I$ and $J$ are maximal ideals.

In this situation, we have $\{0\}\subset I\subset R$ and $\{0\}\subset
J\subset R$, with $I+J=R$. Therefore, $Id\left( R\right) $ is a lattice as
in Figure 1, A\textit{.} We have that $I$ and $J$ are coprime. It is clear
that $I\otimes J=I\cap J=\{0\}$ and $I^{2}=I,J^{2}=J$. From here, it results
that $Ann\left( I\right) =J$ and $Ann\left( J\right) =I$. Also, $%
J\rightarrow I=\left( I:J\right) =I$ and $I\rightarrow J=\left( J:I\right)
=J $. Therefore, we obtain an MV-algebra structure with the following
operations: 
\begin{equation}
\begin{tabular}{l|llll}
$\rightarrow $ & $0$ & $I$ & $J$ & $R$ \\ \hline
$0$ & $R$ & $R$ & $R$ & $R$ \\ 
$I$ & $J$ & $R$ & $J$ & $R$ \\ 
$J$ & $I$ & $I$ & $R$ & $R$ \\ 
$R$ & $0$ & $I$ & $J$ & $R$%
\end{tabular}%
\ \ 
\begin{tabular}{l|llll}
$\otimes $ & $0$ & $I$ & $J$ & $R$ \\ \hline
$0$ & $0$ & $0$ & $0$ & $0$ \\ 
$I$ & $0$ & $I$ & $0$ & $I$ \\ 
$J$ & $0$ & $0$ & $J$ & $J$ \\ 
$R$ & $0$ & $I$ & $J$ & $R$%
\end{tabular}%
\text{.}  \tag{3.3.}
\end{equation}

\textbf{Case ii.} The lattice of ideals is a chain, $Id\left( R\right)
=\{\{0\}\subset I\subset J\subset R\}$, as in Figure 1, \textit{B}.

Since $\{0\}\subset I\subset J\subset R$, the ring is local. Also, we remark
that $I$ is the only minimal ideal. Let $x\in I,x\neq 0$. Therefore, $%
0\not=<x>\subseteq I$. Since $I$ is minimal, we have $I=<x>,$ then $I$ is
finitely generated. We have the following subcases.

1) $I^{2}=\{0\},J^{2}=I,I\otimes J=\{0\}$. We have $Ann\left( I\right)
=\left( 0:I\right) =J$ and $Ann\left( J\right) =\left( 0:J\right) =I$. Also, 
$J\rightarrow I=\left( I:J\right) =J$. Therefore, we obtain an MV-algebra
structure with the following operations:

\begin{equation}
\begin{tabular}{l|llll}
$\rightarrow $ & $0$ & $I$ & $J$ & $R$ \\ \hline
$0$ & $R$ & $R$ & $R$ & $R$ \\ 
$I$ & $J$ & $R$ & $R$ & $R$ \\ 
$J$ & $I$ & $J$ & $R$ & $R$ \\ 
$R$ & $0$ & $I$ & $J$ & $R$%
\end{tabular}%
\ \ 
\begin{tabular}{l|llll}
$\otimes $ & $0$ & $I$ & $J$ & $R$ \\ \hline
$0$ & $0$ & $0$ & $0$ & $0$ \\ 
$I$ & $0$ & $0$ & $0$ & $I$ \\ 
$J$ & $0$ & $0$ & $I$ & $J$ \\ 
$R$ & $0$ & $I$ & $J$ & $R$%
\end{tabular}%
\text{.}  \tag{3.4.}
\end{equation}

2) $I^{2}=\{0\},J^{2}=I,I\otimes J=I$. ~In this case we do not obtaine a
residuated lattice, since $\otimes $ is not associative. For example, $%
I\otimes J^{2}=I^{2}=\{0\}$ and $\left( I\otimes J\right) \otimes J=I\otimes
J=I$.

3) $I^{2}=\{0\},J^{2}=J,I\otimes J=I$. We have $Ann\left( I\right) =\left(
0:I\right) =I$ and $Ann\left( J\right) =\left( 0:J\right) =\{0\}$. Also, $%
J\rightarrow I=\left( I:J\right) =I$. Therefore, for $Id\left( R\right) ~$we
obtain a BL-algebra structure (which is not an MV-algebra) with the
following implication and multiplication tables

\begin{equation}
\begin{tabular}{l|llll}
$\rightarrow $ & $0$ & $I$ & $J$ & $R$ \\ \hline
$0$ & $R$ & $R$ & $R$ & $R$ \\ 
$I$ & $I$ & $R$ & $R$ & $R$ \\ 
$J$ & $0$ & $I$ & $R$ & $R$ \\ 
$R$ & $0$ & $I$ & $J$ & $R$%
\end{tabular}%
\ \ 
\begin{tabular}{l|llll}
$\otimes $ & $0$ & $I$ & $J$ & $R$ \\ \hline
$0$ & $0$ & $0$ & $0$ & $0$ \\ 
$I$ & $0$ & $0$ & $I$ & $I$ \\ 
$J$ & $0$ & $I$ & $J$ & $J$ \\ 
$R$ & $0$ & $I$ & $J$ & $R$%
\end{tabular}%
\text{.}  \tag{3.5.}
\end{equation}

4) $I^{2}=\{0\},J^{2}=J,I\otimes J=\{0\}$. We have $Ann\left( I\right)
=\left( 0:I\right) =J$ and $Ann\left( J\right) =\left( 0:J\right) =I$. Also, 
$J\rightarrow I=\left( I:J\right) =I$. Therefore, we have $J\otimes
(J\rightarrow I)=J\otimes I=\{0\}$ and $J\wedge I=I$, false. Condition $%
\left( div\right) $ is not satisfied. It results that $Id\left( R\right) $
is not a BL-algebra.

5) $I^{2}=\{0\},J^{2}=\{0\},I\otimes J=\{0\}$. We have $Ann\left( I\right)
=\left( 0:I\right) =J$ and $Ann\left( J\right) =\left( 0:J\right) =J$. Also, 
$J\rightarrow I=\left( I:J\right) =J$. Therefore, we have $J\otimes
(J\rightarrow I)=J\otimes J=\{0\}$ and $J\wedge I=I$, false. Condition $%
\left( div\right) $ is not satisfied. It results that $Id\left( R\right) $
is not a BL-algebra.

6) $I^{2}=\{0\},J^{2}=\{0\},I\otimes J=I$, it is not possible, since $%
I=I\otimes J\subset J^{2}=\{0\}$.

7) $I^{2}=I,J^{2}=I,I\otimes J=I$. We have $Ann\left( I\right) =\left(
0:I\right) =0$ and $Ann\left( J\right) =\left( 0:J\right) =0$. Also, $%
J\rightarrow I=\left( I:J\right) =J$. Therefore, for $Id\left( R\right) ~$we
obtain a BL-algebra structure(which is not an MV-algebra) with the following
implication and multiplication tables:

\begin{equation}
\begin{tabular}{l|llll}
$\rightarrow $ & $0$ & $I$ & $J$ & $R$ \\ \hline
$0$ & $R$ & $R$ & $R$ & $R$ \\ 
$I$ & $0$ & $R$ & $R$ & $R$ \\ 
$J$ & $0$ & $J$ & $R$ & $R$ \\ 
$R$ & $0$ & $I$ & $J$ & $R$%
\end{tabular}%
\ \ 
\begin{tabular}{l|llll}
$\otimes $ & $0$ & $I$ & $J$ & $R$ \\ \hline
$0$ & $0$ & $0$ & $0$ & $0$ \\ 
$I$ & $0$ & $I$ & $I$ & $I$ \\ 
$J$ & $0$ & $I$ & $I$ & $J$ \\ 
$R$ & $0$ & $I$ & $J$ & $R$%
\end{tabular}%
\text{.}  \tag{3.6.}
\end{equation}

8) $I^{2}=I,J^{2}=I,I\otimes J=\{0\}$, it is not possible, since $%
I=I^{2}\subset I\otimes J=\{0\}$.

9) $I^{2}=I,J^{2}=J,I\otimes J=I$. We have $Ann\left( I\right) =\left(
0:I\right) =0$ and $Ann\left( J\right) =\left( 0:J\right) =0$. Also, $%
J\rightarrow I=\left( I:J\right) =I$. Therefore, for $Id\left( R\right) ~$we
obtain a BL-algebra structure (which is not an MV-algebra) with the
following implication and multiplication tables:%
\begin{equation}
\begin{tabular}{l|llll}
$\rightarrow $ & $0$ & $I$ & $J$ & $R$ \\ \hline
$0$ & $R$ & $R$ & $R$ & $R$ \\ 
$I$ & $0$ & $R$ & $R$ & $R$ \\ 
$J$ & $0$ & $I$ & $R$ & $R$ \\ 
$R$ & $0$ & $I$ & $J$ & $R$%
\end{tabular}%
\ \ 
\begin{tabular}{l|llll}
$\otimes $ & $0$ & $I$ & $J$ & $R$ \\ \hline
$0$ & $0$ & $0$ & $0$ & $0$ \\ 
$I$ & $0$ & $I$ & $I$ & $I$ \\ 
$J$ & $0$ & $I$ & $J$ & $J$ \\ 
$R$ & $0$ & $I$ & $J$ & $R$%
\end{tabular}%
\text{.}  \tag{3.7.}
\end{equation}

10) $I^{2}=I,J^{2}=J,I\otimes J=\{0\}$, it is not possible, since $%
I=I^{2}\subset I\otimes J=\{0\}$.\medskip

11) $I^{2}=I,J^{2}=\{0\},I\otimes J=I$, it is not possible, since $%
I=I\otimes J\subset J^{2}=\{0\}$.\medskip

12) $I^{2}=I,J^{2}=\{0\},I\otimes J=\{0\}$, it is not possible, since from $%
I\subset J$ we obtain $I^{2}\subset J^{2}=\{0\}$, so $I=\{0\},$ a
contradiction.\medskip

\textbf{Remark 3.7.} Counting the BL-algebras of order $4$ from the above
proposition, we obtain $5$ BL-algebras of order $4$, two MV and three
BL-chain. In this way, we recover the results obtained in [FP; 22-2], Table
2. From the above remark, we see that all BL-algebras of order four can be
obtained as a lattice of ideals of a commutative unitary ring $R$. For
MV-algebra given by the relation $\left( 3.3\right) $, there is an example
even of a ring which this MV-algebras are isomorphic to, namely $R=\mathbb{Z}%
_{2}\times \mathbb{Z}_{2}$, and for MV-algebra given by the relation $\left(
3.4\right) $ an example can be the ring $\mathbb{Z}_{8\text{ \ }}$(see [FP;
22-1], Table 2). It is interesting to search and find, if there are,
examples of commutative unitary rings $R$ with $Id\left( R\right) $ a
BL-algebra given by relations $\left( 3.5\right) ,\left( 3.6\right) ,\left(
3.7\right) $.\medskip

\textbf{Remark 3.8.} a) Other examples of BL-rings are presented in the
following. From ([HLN; 18], Example 2.4 and Remark 2.5) we know that a
Noetherian multiplication ring is a BL-ring. The ring $\mathbb{\,Z}$ is a
Noetherian multiplication ring ( it is principal). Therefore, from the above
properties $(\mathbb{Z}_{n},+,\cdot )$ is a Noetherian multiplication ring.

b) We consider the ring $(\mathbb{Z}_{n},+,\cdot )$ with $%
n=p_{1}p_{2}...p_{r},p_{1},p_{2},...,p_{r}$ being distinct prime numbers, $%
r\geq 2$. The ring $R=\mathbb{Z}_{n}[X]/\left( f\right) ,$ with $f$ a
polynomial mod $n$ of degree $q\geq 2$, is a Noetherian multiplication ring,
therefore a BL-ring. From here, we have that the algebra of ideals, $%
Id\left( R\right) $, are BL-algebra which is not a BL-chain ($\mathbb{Z}_{n}$
is a direct product of the fields $\mathbb{Z}_{p_{i}}$).

c) If $n=p^{r},$ $r\geq 2$, there is an MV-algebra isomorphic to $\mathbb{Z}%
_{n}$ ([FP; 22-1]).

d) From ([AB; 19], Corollary 1.3), we know that an Artinian ring is a
multiplication ring if and only if it is a finite product of Artinian local
principal ideal rings. Moreover, if $r=1,$the ring $R=\mathbb{Z}%
_{p}[X]/\left( f\right) $ is Artinian (is finite), local and principal ring,
therefore it is a BL-ring.

e) If we take the ring $\mathbb{Z}_{p}^{r}=$ $\underset{r-time}{\underbrace{%
\mathbb{Z}_{p}\times ...\times \mathbb{Z}_{p}}}$, we have that the ring $R=%
\mathbb{Z}_{p}^{r}[X]/\left( f\right) $ is a Noetherian multiplication ring,
therefore a BL-ring.\medskip

\textbf{Remark 3.9.} Let $R$ be a commutative unitary ring which has exactly
five ideals. Therefore, the lattice $Id\left( R\right) $ can be of the form $%
A$-$E$ from Figure 2 and Figure 3. Since $5$ is a prime number, with the
above notations, we have only commutative unitay rings $R$ of the form $%
\left( n_{m}\left( R\right) ,n_{p}\left( R\right) ,n_{I}\left( R\right)
\right) =\left( 1,1,5\right) $ with $Id\left( R\right) $ a BL-algebra, that
means $Id\left( R\right) $ is a chain (Figure E) or has Hasse diagram as in
Figure C. Also, we remark that a BL-algebra cannot have the form as in
Figure C since, in this case $(a\rightarrow b)\vee (b\rightarrow a)=c\neq 1,$
so (prel) condition is not satisfied. We conclude that $Id\left( R\right) $
is a chain.

In [NL; 03], Corollary 28, the authors proved that each finite BL-algebras
are isomorphic to a direct product of BL-comets. A BL-comet of order $5$ has
Hasse diagram as in Figures 3, \textit{D or E}. If we consider a BL-comet of
order 5 as in Figure 3, D, generated by a commutative ring $R,$ then $%
Id\left( R\right) $ has two maximal ideals, impossible. We conclude that,
BL-comets $Id\left( R\right) $ generated by commutative rings $R$ with five
ideals are chains. As in Remark 3.6, we can obtain all $8$ BL-chain of order 
$5$, one MV and seven BL-chain as $Id\left( R\right) $ of a ring with $5$
ideals with $\left( n_{m}\left( R\right) ,n_{p}\left( R\right) ,n_{I}\left(
R\right) \right) =\left( 1,1,5\right) $. In this way, we recover the results
obtained in [FP; 22-2], Table 2.

\begin{figure}[tbph]
\centerline{\includegraphics[width=3in, height=1.3in]{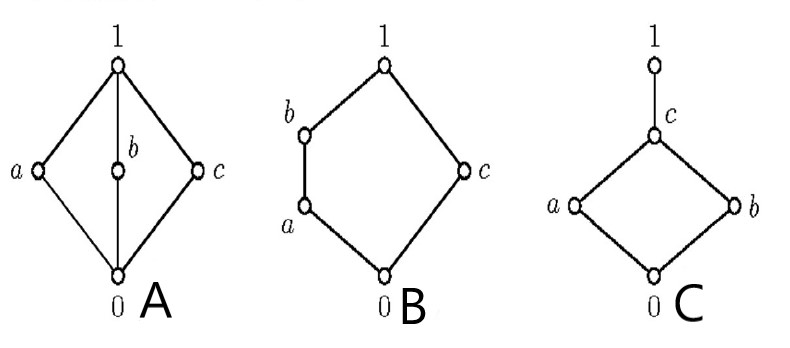}}
\caption{Latices with five elements.}
\end{figure}
\begin{figure}[tbph]
\centerline{\includegraphics[width=2.8in, height=1.3in]{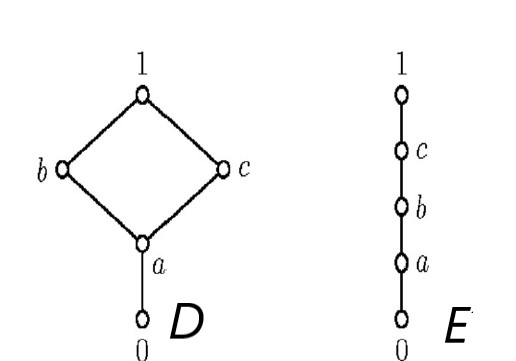}}
\caption{Latices with five elements.}
\end{figure}

\textbf{Remark 3.10.} a) With the above notations, we do not have examples
of commutative unitary rings $R$ of the form $\left( n_{m}\left( R\right)
,n_{p}\left( R\right) ,n_{I}\left( R\right) \right) =\left( 2,2,6\right) $
with $Id\left( R\right) $ a BL-algebra which is not an MV-algebra. Indeed,
if such a ring exists, it is an Artinian ring, therefore it is a product of
Artinian local nonzero rings, $R=R_{1}\times R_{2}$. One of these rings, $%
R_{1},$ has two ideals, and the other one, $R_{2},$ has three ideals.
Therefore $Id\left( R_{1}\right) $ and $Id\left( R_{2}\right) $ must be
BL-algebras and at least one should not be an MV-algebra.. But, the only
BL-algebra with 2 elements, $Id(Z_{2}),$ is an MV-algebra and from
Proposition 3.1, we have that there are not commutative unitary rings with
three ideals such that $Id\left( R_{2}\right) $ is a BL-algebra which are
not an MV-algebra. In this case we have only finite commutative unitay rings 
$R$ of the form $\left( n_{m}\left( R\right) ,n_{p}\left( R\right)
,n_{I}\left( R\right) \right) =\left( 1,1,6\right) $ with $Id\left( R\right) 
$ a BL-algebra which is not an MV-algebra.

b) We can have only commutative unitay $R$ rings of the form $\left(
n_{m}\left( R\right) ,n_{p}\left( R\right) ,n_{I}\left( R\right) \right)
=\left( 1,1,7\right) $ with $Id\left( R\right) $ a BL-algebra, since $7$ is
a prime number, that means $Id\left( R\right) $ is a chain.

b) If $n_{I}\left( R\right) =8$, we can have commutative unitary rings of
the form $\left( n_{m}\left( R\right) ,n_{p}\left( R\right) ,n_{I}\left(
R\right) \right) =\left( 2,2,8\right) \,\ $or $\left( n_{m}\left( R\right)
,n_{p}\left( R\right) ,n_{I}\left( R\right) \right) =\left( 3,3,8\right) $,
with $Id\left( R\right) $ a BL-algebra.

d) For $n_{I}\left( R\right) =9,$ we do not have commutative unitary rings
of the form $\left( n_{m}\left( R\right) ,n_{p}\left( R\right) ,n_{I}\left(
R\right) \right) =\left( 2,2,9\right) $, with $Id\left( R\right) $ a
BL-algebra which is not an MV-algebra, since there are not commutative
unitary rings $R$ with three ideals such that $Id\left( R\right) $ is a
BL-algebra which is not an MV-algebra. But we can have rings $R$ of the form 
$\left( n_{m}\left( R\right) ,n_{p}\left( R\right) ,n_{I}\left( R\right)
\right) =\left( 1,1,9\right) $, that means $Id\left( R\right) $ is a chain.
For MV-algebras of order $9$, there are examples of rings whit which this
MV-algebras is isomorphic to, namely $\mathbb{Z}_{3}\times \mathbb{Z}_{3}$
or $\mathbb{Z}_{9}$(see [FP; 22-1], Table 2).

e) For $n_{I}\left( R\right) =10,$ we can have $\left( n_{m}\left( R\right)
,n_{p}\left( R\right) ,n_{I}\left( R\right) \right) =\left( 2,2,10\right) $
or $\left( n_{m}\left( R\right) ,n_{p}\left( R\right) ,n_{I}\left( R\right)
\right) =\left( 1,1,10\right) $ such that $Id\left( R\right) $ is a
BL-algebra which is not an MV-algebra.\medskip

\textbf{Example 3.11. }Let $R$ be a BL ring with $Id\left( R\right) $
infinite and $I$ be an idempotent ideal in $R,I^{2}=I$. From Example 3.5, we
know that such a ring exists and $\{0,I,R\}$ forms a BL-subalgebra of order $%
3$ with the multiplication tables given by the relation $\left( 3.2\right) $%
. We consider the ring $R\times R$, therefore 
\begin{eqnarray*}
\mathcal{B} &=&\{\left( 0,0\right) ,\left( 0,I\right) ,\left( 0,R\right)
,\left( I,0\right) ,\left( I,I\right) ,\left( I,R\right) ,\left( R,0\right)
,\left( R,I\right) ,\left( R,R\right) \} \\
\newline
&=&\{O,A,B,C,D,E,F,G,Z\}
\end{eqnarray*}%
$~$is a BL-subalgebra of order $9$ (see Figure 4). We have the following
implication and multiplication tables:

\begin{equation}
\begin{tabular}{l|lllllllll}
$\rightarrow $ & $O$ & $A$ & $B$ & $C$ & $D$ & $E$ & $F$ & $G$ & $Z$ \\ 
\hline
$O$ & $Z$ & $Z$ & $Z$ & $Z$ & $Z$ & $Z$ & $Z$ & $Z$ & $Z$ \\ 
$A$ & $F$ & $Z$ & $Z$ & $F$ & $Z$ & $Z$ & $F$ & $Z$ & $Z$ \\ 
$B$ & $F$ & $G$ & $Z$ & $F$ & $G$ & $Z$ & $F$ & $G$ & $Z$ \\ 
$C$ & $B$ & $B$ & $B$ & $\mathbf{Z}$ & $\mathbf{Z}$ & $\mathbf{Z}$ & $%
\mathbf{Z}$ & $\mathbf{Z}$ & $\mathbf{Z}$ \\ 
$D$ & $O$ & $B$ & $B$ & $\mathbf{F}$ & $\mathbf{Z}$ & $\mathbf{Z}$ & $%
\mathbf{F}$ & $\mathbf{Z}$ & $\mathbf{Z}$ \\ 
$E$ & $O$ & $A$ & $B$ & $\mathbf{F}$ & $\mathbf{G}$ & $\mathbf{Z}$ & $%
\mathbf{F}$ & $\mathbf{G}$ & $\mathbf{Z}$ \\ 
$F$ & $B$ & $B$ & $B$ & $\mathbf{E}$ & $\mathbf{E}$ & $\mathbf{E}$ & $%
\mathbf{Z}$ & $\mathbf{Z}$ & $\mathbf{Z}$ \\ 
$G$ & $O$ & $B$ & $B$ & $\mathbf{C}$ & $\mathbf{E}$ & $\mathbf{E}$ & $%
\mathbf{F}$ & $\mathbf{Z}$ & $\mathbf{Z}$ \\ 
$Z$ & $O$ & $A$ & $B$ & $\mathbf{C}$ & $\mathbf{D}$ & $\mathbf{E}$ & $%
\mathbf{F}$ & $\mathbf{G}$ & $\mathbf{Z}$%
\end{tabular}%
\ \ \ \ \ 
\begin{tabular}{l|lllllllll}
$\otimes $ & $O$ & $A$ & $B$ & $C$ & $D$ & $E$ & $F$ & $G$ & $Z$ \\ \hline
$O$ & $O$ & $O$ & $O$ & $O$ & $O$ & $O$ & $O$ & $O$ & $O$ \\ 
$A$ & $O$ & $A$ & $A$ & $O$ & $A$ & $A$ & $O$ & $A$ & $A$ \\ 
$B$ & $O$ & $A$ & $B$ & $O$ & $A$ & $B$ & $O$ & $A$ & $B$ \\ 
$C$ & $O$ & $O$ & $O$ & $\mathbf{C}$ & $\mathbf{C}$ & $\mathbf{C}$ & $%
\mathbf{C}$ & $\mathbf{C}$ & $\mathbf{C}$ \\ 
$D$ & $O$ & $A$ & $A$ & $\mathbf{C}$ & $\mathbf{D}$ & $\mathbf{D}$ & $%
\mathbf{C}$ & $\mathbf{D}$ & $\mathbf{D}$ \\ 
$E$ & $O$ & $A$ & $B$ & $\mathbf{C}$ & $\mathbf{D}$ & $\mathbf{E}$ & $%
\mathbf{C}$ & $\mathbf{D}$ & $\mathbf{E}$ \\ 
$F$ & $O$ & $O$ & $O$ & $\mathbf{C}$ & $\mathbf{C}$ & $\mathbf{C}$ & $%
\mathbf{F}$ & $\mathbf{F}$ & $\mathbf{F}$ \\ 
$G$ & $O$ & $A$ & $A$ & $\mathbf{C}$ & $\mathbf{D}$ & $\mathbf{D}$ & $%
\mathbf{F}$ & $\mathbf{G}$ & $\mathbf{G}$ \\ 
$Z$ & $O$ & $A$ & $B$ & $\mathbf{C}$ & $\mathbf{D}$ & $\mathbf{E}$ & $%
\mathbf{F}$ & $\mathbf{G}$ & $\mathbf{Z}$%
\end{tabular}
\tag{3.8.}
\end{equation}

\begin{figure}[tbph]
\centerline{\includegraphics[width=2in, height=1.3in]{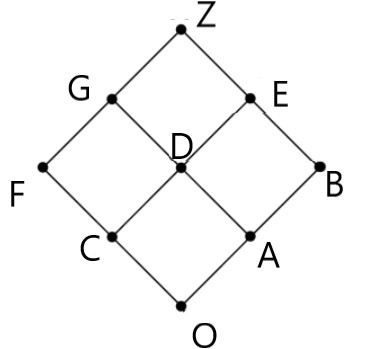}}
\caption{Latices with nine elements.}
\end{figure}
From here, we remark that $\{C,D,E,F,G,Z\}$ is a BL-algebra of order $6~$(
see bold parts from relation $\left( 3.8\right) $), extracted from the
algebra $\mathcal{B}$, $\{O,D,E,G,Z\}$ is a BL-subalgebra of order $5$ (it
is a BL-comet, from the below tables) and $\{F,G,Z\}$ is a BL-algebra of
order $3~$(as in relation $\left( 3.7\right) $), extracted from the algebra $%
\mathcal{B}$ (see Figure 5).\ 
\begin{equation*}
\begin{tabular}{l|lllll}
$\rightarrow $ & $O$ & $D$ & $E$ & $G$ & $Z$ \\ \hline
$O$ & $Z$ & $Z$ & $Z$ & $Z$ & $Z$ \\ 
$D$ & $O$ & $Z$ & $Z$ & $Z$ & $Z$ \\ 
$E$ & $O$ & $G$ & $Z$ & $G$ & $Z$ \\ 
$G$ & $O$ & $E$ & $E$ & $Z$ & $Z$ \\ 
$Z$ & $O$ & $D$ & $E$ & $G$ & $Z$%
\end{tabular}%
\ \ \ 
\begin{tabular}{l|lllll}
$\otimes $ & $O$ & $D$ & $E$ & $G$ & $Z$ \\ \hline
$O$ & $O$ & $O$ & $O$ & $O$ & $O$ \\ 
$D$ & $O$ & $D$ & $D$ & $D$ & $D$ \\ 
$E$ & $O$ & $D$ & $E$ & $D$ & $E$ \\ 
$G$ & $O$ & $D$ & $D$ & $G$ & $G$ \\ 
$Z$ & $O$ & $D$ & $E$ & $G$ & $Z$%
\end{tabular}%
\text{.}
\end{equation*}

\begin{figure}[tbph]
\centerline{\includegraphics[width=2in, height=1.3in]{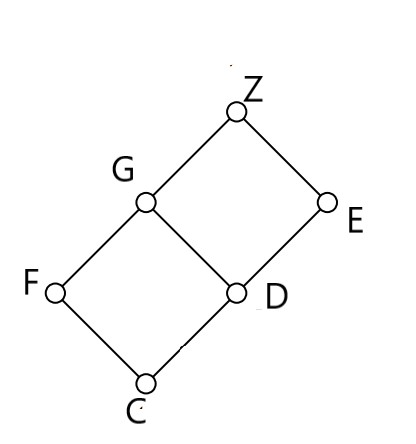}}
\caption{Latices with six elements.}
\end{figure}
\textbf{Conclusions.} In this paper, we presented some rings with a finite
number of ideals and their connections with BL-rings. BL-rings are
commutative unitary rings whose lattice of ideals can be equipped with a
structure of BL-algebra. We proved that for a given number of ideals, some
of these rings $R$ such that $Id\left( R\right) $ is a BL-algebra can exist.
As a further research we intend to find more examples of such type of rings $%
R$ with $Id\left( R\right) $ a finite BL-algebra.

\begin{equation*}
\end{equation*}

\textbf{References}%
\begin{equation*}
\end{equation*}

[AB; 19] Alsuraiheed, T., Bavula, V. V., \textit{Characterization of
multiplication commutative rings with finitely many minimal prime ideals},
Communications in Algebra, 47(11))(2019), 4533-4540.

[A; 76] Anderson D. D., \textit{Multiplication ideals, multiplication rings,
and the ring} $R(X)$, Can. J. Math., Vol. 28(4)(1976), 760-768.

[AF; 92] Anderson, F. W., Fuller, K., (1992), \textit{Rings and categories
of modules}, Graduate Texts in Mathematics, 13(1992), 2 ed.,
Springer-Verlag, New York.

[AM; 69] Atiyah, M. F., MacDonald, I. G., I\textit{ntroduction to
Commutative Algebra}, Addison-Wesley Publishing Company, London, 1969.

[CHA; 58] Chang, C.C.,\textit{\ Algebraic analysis of many-valued logic},
Trans. Amer. Math. Soc. 88(1958), 467-490.

[FK; 12] Filipowicz, M., Kepczyk, M., \textit{A note on zero-divisors of
commutative rings}, Arab J Math, 1(2012), 191--194.

[FP; 22-1] Flaut, C., Piciu, D., \textit{Connections between commutative
rings and some algebras of logic}, Iranian Journal of Fuzzy Systems,
19(6)(2022), 93-110.

[FP; 22-2] Flaut, C., Piciu, D., \textit{Some Examples of BL-Algebras Using
Commutative Rings}, Mathematics 10(2022), 4739, 1-15.

[HLN; 18] Heubo-Kwegna, O. A., Lele, C., Ndjeya, S., Nganou, J. B., \textit{%
BL-rings}, Logic Journal of the IGPL, 26(3)(2018), 290--299.

[H; 98] H\'{a}jek, P., \textit{Metamathematics of Fuzzy Logic}, Trends in
Logic-Studia Logica Library 4, Dordrecht: Kluwer Academic Publishers \textbf{%
1998.\medskip }

[NL; 03] Di Nola, A., Lettieri, A., \textit{Finite BL-algebras}, Discrete
Mathematics, 269(2003), 93-112.

[NL; 05] Di Nola, A., Lettieri, A., \textit{Finiteness based results in
BL-algebras}, Soft Comput 9(2005), 889--896.

[P; 07] Piciu, D., Algebras of fuzzy logic, Ed. Universitaria, Craiova, 2007.

[TT;22] Tchoffo Foka, S. V., Tonga, M., \textit{Rings and residuated
lattices whose fuzzy ideals form a Boolean algebra}, Soft Computing, 26
(2022) 535-539.

[T; 99] Turunen, E., \textit{Mathematics Behind Fuzzy Logic},
Physica-Verlag, \textbf{1999}.

[WD; 39] Ward, M., Dilworth, R.P., \textit{Residuated lattices}, Trans. Am.
Math. Soc. 45(1939), 335--354.

\bigskip

\begin{equation*}
\end{equation*}

{\small Mariana Floricica C\u{a}lin}

{\small Faculty of Psychology and Educational Science, }

{\small Ovidius University of Constan\c{t}a, Rom\^{a}nia,}

{\small Bd. Mamaia 124, 900527,}

{\small e-mail: fmarianacalin@gmail.com\medskip }

\bigskip

Cristina Flaut

{\small Faculty of Mathematics and Computer Science, Ovidius University,}

{\small Bd. Mamaia 124, 900527, Constan\c{t}a, Rom\^{a}nia,}

{\small \ http://www.univ-ovidius.ro/math/}

{\small e-mail: cflaut@univ-ovidius.ro; cristina\_flaut@yahoo.com}

\bigskip

Dana Piciu

{\small Faculty of \ Science, University of Craiova, }

{\small A.I. Cuza Street, 13, 200585, Craiova, Romania,}

{\small http://www.math.ucv.ro/dep\_mate/}

{\small e-mail: dana.piciu@edu.ucv.ro, piciudanamarina@yahoo.com}

\end{document}